\documentclass[a4paper, 12pt]{amsart}
\usepackage{setspace}
\usepackage{amsmath,amssymb,amsthm}
\usepackage{slashbox}
\usepackage{float}
\usepackage{hhline}
\usepackage{mathrsfs}
\usepackage[capitalize]{cleveref}
\usepackage{autonum}

\newtheorem{theorem}{Theorem}[section]
\crefname{theorem}{Theorem}{Theorems}
\newtheorem{corollary}[theorem]{Corollary}
\crefname{corollary}{Corollary}{Corollaries}

\crefname{fact}{Fact}{Facts}

\crefname{proposition}{Proposition}{Propositions}

\crefname{lemma}{Lemma}{Lemmas}

\crefname{conjecture}{Conjecture}{Conjectures}

\crefname{question}{Question}{Questions}
\theoremstyle{definition}

\crefname{definition}{Definition}{Definitions}

\crefname{example}{Example}{Examples}
\newtheorem{remark}[theorem]{Remark}
\crefname{remark}{Remark}{Remarks}
\newcommand{\R}{{\mathbb R}}
\newcommand{\Z}{{\mathbb Z}}
\newcommand{\N}{{\mathbb N}}
\newcommand{\C}{{\mathbb C}}

\newcommand{\sm}{\setminus}
\newcommand{\I}{\sqrt{-1}}
\newcommand{\wt}{{\rm wt}}

\def\ip<#1>{\langle #1\rangle}
\newcommand{\suma}[2]{\sideset{}{^*}\sum_{#1}^{#2}}

\allowdisplaybreaks[1]
\makeatletter

\@addtoreset{equation}{section}
\makeatother

\title{On the parity result for multiple Dirichlet series}
\author[S.~KADOTA]{Shin-ya Kadota}
\address{National Institute of Technology, Niihama College, Yagumocho, Ehime 792-8580, Japan.}
\email{kadota@sci.niihama-nct.ac.jp}
\date{}
\begin{document}
\maketitle
\begin{abstract}
In this paper, we discuss the parity result for multiple Dirichlet series which contains some special values of multiple zeta functions as special cases, such as Mordell--Tornheim type of multiple zeta values, zeta values of the root systems and so on. Moreover, we can give an explicit expression in terms of lower series by using the main theorem.
\end{abstract}
\section{Introduction}
Parity results are one of the important properties of special values of multiple zeta functions, and have been extensively researched. For example, the parity result for Euler--Zagier type of multiple zeta values was proved by Goncharov \cite[Corollary 7.2]{Goncharov98}, Tsumura \cite{Tsumura04AA} and Ihara, Kaneko and Zagier \cite{IKZ06}. For an $r$-tuple of positive integers ${\bf k}=(k_1, k_2, \dots, k_r)~(k_r\ge2)$, the Euler--Zagier type of multiple zeta value is defined by
\begin{align}
\zeta_{r, EZ}({\bf k})=\sum_{m_1=1}^\infty\dots\sum_{m_r=1}^\infty\frac{1}{m_1^{k_1}(m_1+m_2)^{k_2}\dots(m_1+m_2+\dots+m_r)^{k_r}}.
\end{align}
For a tuple of numbers ${\bf k}=(k_1, k_2, \dots, k_r)$, we define $\wt({\bf k})$ as the sum of all entries of ${\bf k}$, namely $\wt({\bf k}):=\sum_{j=1}^rk_j$ which is called the weight. The parity result for Euler--Zagier type of multiple zeta values is the following property.
\begin{theorem}
When $\wt({\bf k})$ and $r\ge2$ are of different parity, $\zeta_{r, EZ}({\bf k})$ can be written as a rational linear combination of the product of Euler--Zagier type of multiple zeta values $\zeta_{s, EZ}({\boldsymbol \ell})~(s<r, \wt({\boldsymbol \ell})<\wt({\bf k}))$.
\end{theorem}
The parity result for Mordell--Tornheim type of multiple zeta values was proved by Tsumura \cite{Tsumura05} and later by Onodera \cite[Theorem 3]{Onodera11} (see also the paper of Bradley and Zhou \cite{BZ10}). The Mordell--Tornheim type of multiple zeta value is defined for an $(r+1)$-tuple of positive numbers ${\bf k}=(k_1, \dots, k_r, k_{r+1})$ by
\begin{align}
\zeta_{r, MT}({\bf k})=\sum_{m_1=1}^\infty\dots\sum_{m_r=1}^\infty\frac{1}{m_1^{k_1}m_2^{k_2}\dots m_r^{k_r}(m_1+m_2+\dots+m_r)^{k_{r+1}}}.
\end{align}
This series with $k_i=s_i\in\C$ was introduced by Matsumoto in \cite{Matsumoto03BMS} first.
\begin{theorem}
When $\wt({\bf k})$ and $r\ge2$ are of different parity, $\zeta_{r, MT}({\bf k})$ can be written as a rational linear combination of the product of Mordell--Tornheim type of multiple zeta values $\zeta_{s, MT}({\boldsymbol \ell})~(s<r, \wt({\boldsymbol \ell})<\wt({\bf k}))$.
\end{theorem}
From these two theorems, we can regard the parity result as a property that $r$-ple series can be written in terms of $s(<r)$-ple series. In that sense, some results on such a property for other objects (for example, zeta values of root systems and multiple polylogarithms) are already known (see \cite{KMT14FACM, KOT17, Nakamura12, Okamoto12, Panzer17, Tsumura04AM}). We introduce the main object of the present paper in the next section.
\section{Preliminaries and main results}
Let $\N$ be the set of positive integers, $\N_0$ be the set of non negative integers, $\Z$ be the set of integers, $\R$ be the set of real numbers, $\C$ be the set of complex numbers. For a positive integer $m$, we define $[m]:=\{1, 2, \dots, m\}.$ Let $r$ and $\ell$ be positive integers and $M(\ell, r, \N_0)$ be the set of $\ell\times r$ matrices $A=(a_{ij})_{\substack{i\in[\ell]\\j\in[r]}}$ whose entries are non negative integers which satisfies the following two conditions.
\begin{enumerate}
\item $A$ has no zero row vectors.
\item $A$ has no zero column vectors.
\end{enumerate}
The main object of the present paper is the following multiple Dirichlet series which is defined for ${\bf h}=(h_1, \dots, h_r)\in\N^r$, ${\bf k}=(k_1, \dots, k_\ell)\in\N^\ell$, ${\bf y}=(y_1, \dots, y_r)\in\R^r$ and $A\in M(\ell, r, \N_0)$ by
\begin{align}
&\zeta_{r, \ell}({\bf h}, {\bf k}, {\bf y}, A)\\
&:=\sum_{m_1=1}^\infty\dots\sum_{m_r=1}^\infty\frac{e^{2\pi \I(m_1y_1+\dots+m_ry_r)}}{m_1^{h_1}\dots m_r^{h_r}(a_{11}m_1+\dots+a_{1r}m_r)^{k_1}\dots(a_{\ell1}m_1+\dots+a_{\ell r}m_r)^{k_\ell}}\\
&=\sum_{m_1=1}^\infty\dots\sum_{m_r=1}^\infty\prod_{j\in[r]}\frac{e(m_jy_j)}{m_j^{h_j}}\prod_{i\in[\ell]}\frac{1}{(a_{i1}m_1+\dots+a_{ir}m_r)^{k_i}},
\end{align}
where $e(\alpha):=\exp(2\pi\sqrt{-1}\alpha)$. The first condition of $A$ guarantees that the denominator does not vanish, and the second condition guarantees that $\zeta_{r, \ell}({\bf h}, {\bf k}, {\bf y}, A)$ converges absolutely.

Here, we introduce notations. For a non empty subset $J=\{j_1<j_2<\dots<j_{|J|}\}\subset[r]$, we define $\bar{J}:=[r]\sm J$. For a tuple of numbers ${\bf h}=(h_1, \dots, h_r)$ and a matrix $A\in M(\ell, r, \N_0)$, we define ${\bf h}_J$ and a sub matrix $A_J$ as follows.
\begin{align}
{\bf h}_J:=(h_j)_{j\in J}=(h_{j_1}, h_{j_2}, \dots, h_{j_m}),\quad A_J:=(a_{ij})_{\substack{i\in[\ell]\\j\in J}}.
\end{align}
For the empty set, we define $\wt({\bf h}_\emptyset)=0$. Moreover, we define the set of numbers of non zero row vectors of $A_J$.
\begin{align}
I=I_J:=\{i\in[\ell]\mid a_{ij}\neq0\mbox{ for some }j\in J \},\quad \bar{I}=\bar{I}_J:=[\ell]\sm I_J.
\end{align}
\par
In order to state the main theorem, we need more preliminaries. Let $m$ be a positive integer. We assume that the real vector space $\R^m$ is equipped with the normal inner product $\ip<\cdot, \cdot>$. We denote $\vec{f}$ the part of $\R^m$ in $f\in\R^m\times\C$ and $\dot{f}$ the part in $\C$ of $f\in\R^m\times\C$. Namely, we express $f\in\R^m\times\C$ as $f=(\vec{f}, \dot{f})$. We define $\ip<V>:=\sum_{{\bf v}\in V}\Z{\bf v}$ for $V\subset\R^m$, and $\vec{W}:=\{\vec{f}\mid f=(\vec{f}, \dot{f})\in W\}$ for $W\subset\R^m\times\C$. Let $\Lambda$ be a subset of $(\Z^m\sm\{\vec{0}\})\times\C$ with satisfying $|\Lambda|<\infty$ and ${\rm rank}\ip<\vec{\Lambda}>=m$. 
For $H\subset\Lambda$ with ${\rm rank}\ip<\vec{H}>=m-1$, we define $\mathfrak{H}_H:=\sum_{g\in H}\R\vec{g}$. For $\Lambda$, let $\mathscr{B}(\Lambda)$ be the set of all subsets $B\subset\Lambda$ such that $\vec{B}$ forms a basis of $\R^m$. For $B=\{f_1, \dots, f_m\}\in\mathscr{B}(\Lambda)$, let $\vec{B}^*=\{\vec{f}_1^B, \dots, \vec{f}_m^B\}$ be the dual basis of $\vec{B}=\{\vec{f}_1, \dots, \vec{f}_m\}$ in $\R^m$.\par
For real numbers, we define the multi dimensional generalization of fractional part $\{\cdot\}$ which was introduced in \cite{KMT10PLMS} by Komori, Matsumoto and Tsumura. Let $\mathscr{R}(\Lambda)$ be the set of all subsets $R=\{g_1, \dots, g_{m-1}\}\subset\Lambda$ satisfying that $\vec{R}=\{\vec{g}_1, \dots, \vec{g}_{m-1}\}$ is a linearly independent set. Here, we fix one vector
\begin{align}
\rho\in\R^m\sm \bigcup_{R\in\mathscr{R}(\Lambda)}\mathfrak{H}_R
\end{align}
satisfying that $\ip<\rho, \vec{f}^B>\neq0$ for all $B\in\mathscr{B}(\Lambda)$ and $f\in B$. According to $\rho$, for ${\bf y}\in\R^m, B\in\mathscr{B}(\Lambda)$ and $f\in B$, we define the multi dimensional fractional part as follows.
\begin{align}
\{{\bf y}\}_{B, f}:=\begin{cases}\{\ip<{\bf y}, \vec{f}^B>\}&{\rm if}~\ip<\rho, \vec{f}^B>>0,\\1-\{-\ip<{\bf y}, \vec{f}^B>\}&{\rm if}~\ip<\rho, \vec{f}^B><0.\end{cases}
\end{align}
\par
Now we define the generating function for a non empty subset $J\subset[r]$, a matrix $A\in M(\ell ,r, \N_0)$ and ${\mathbf m}_{\bar{J}}\in\N^{|\bar{J}|}$. We put
\begin{align}
f_j&=f_{J, j}=((\delta_{jk})_{k\in J}, 0)\in\R^{|J|}\times\C\quad(j\in J),\\
f_{r+i}&=f_{J, r+i}=((a_{ij})_{j\in J}, -\sum_{j\in\bar{J}}a_{ij}m_j)\in\R^{|J|}\times\C\quad(i\in [\ell]),
\end{align}
and we define $f({\bf m}_J):=\ip<\vec{f}, {\bf m}_J>+\dot{f}$ for the above vectors $f$ and ${\bf m}_J=(m_j)_{j\in J}\in\N^{|J|}$, where $\delta_{jk}$ denotes Kronecker's delta symbol.
\begin{remark}\label{remark:1}
We can rewrite the set $I$ as $\{i\in[\ell]\mid\vec{f}_{r+i}\neq{\bf 0}\}$.
\end{remark}
Moreover, we put
\begin{align}
\Lambda=\Lambda_J=\{f_j, f_{r+i}\mid j\in J, i\in I\}.
\end{align}
It is easy to see that ${\rm rank}\ip<\vec{\Lambda}>=|J|$. For ${\bf y}_J\in\R^{|J|}$ and ${\bf t}_{\Lambda}=(t_f)_{f\in\Lambda}\in\C^{|\Lambda|}$, we define
\begin{align}
G({\bf t}_{\Lambda}, {\bf y}_J; \Lambda)&=\sum_{B\in\mathscr{B}(\Lambda)}\left(\prod_{g\in\Lambda\sm B}\frac{-t_g}{\dot{g}-\sum_{f\in B}\dot{f}\ip<\vec{g}, \vec{f}^B>-(t_g-\sum_{f\in B}t_f\ip<\vec{g}, \vec{f}^B>)}\right)\\
&\quad\times\frac{1}{|\Z^{|J|}/\ip<\vec{B}>|}\sum_{{\bf w}\in \Z^{|J|}/\ip<\vec{B}>}\left(\prod_{f\in B}\frac{2\pi\I t_fe((t_f-\dot{f})\{{\bf y}_J+{\bf w}\}_{B, f})}{e(t_f)-1}\right),
\end{align}
and denote by $D({\bf h}_J, {\bf k}_I, {\bf y}_J; \Lambda)$ the coefficients of the Taylor expansion of $G({\bf t}_{\Lambda}, {\bf y}_J; \Lambda)$ around the origin in ${\bf t}_{\Lambda}$. Namely,
\begin{align}
G({\bf t}_{\Lambda}, {\bf y}_J; \Lambda)=\sum_{\substack{{\bf h}_J\in\N_0^{|J|}\\{\bf k}_{I}\in\N_0^{|I|}}}D({\bf h}_J, {\bf k}_I, {\bf y}_J; \Lambda)\prod_{\substack{j\in J\\i\in I}}\frac{t_j^{h_j}t_{r+i}^{k_i}}{h_j!k_i!},
\end{align}
where $t_f$ denotes $t_k$ for $f=f_k\in\Lambda$. Now, we state main theorem.

\begin{theorem}\label{theorem:main theorem parity}
For $A\in M(\ell, r; \N_0), {\bf h}=(h_1, \dots, h_r)\in\N^r, {\bf k}=(k_1, \dots, k_\ell)\in\N^\ell$ satisfying that
\begin{align}
\mathop{\sum_{m_1=1}^\infty\dots\sum_{m_r=1}^\infty}_{\sum_{j\in J}a_{ij}m_j-\sum_{j\in\bar{J}}a_{ij}m_j\neq0~(i\in[\ell])}\prod_{j\in [r]}\frac{1}{m_j^{h_j}}\prod_{i\in[\ell]}\frac{1}{|\sum_{j\in J}a_{ij}m_j-\sum_{j\in\bar{J}}a_{ij}m_j|^{k_i}}
\end{align}
converges for any subset $J\subset[r]$, and any ${\bf y}\in\R^r$, we have
\begin{align}
\zeta_{r, \ell}({\bf h}, {\bf k}, {\bf y}, A)+(-1)^{\wt({\bf h})+\wt({\bf k})+r+1}\zeta_{r, \ell}({\bf h}, {\bf k}, -{\bf y}, A)=\sum_{\emptyset\neq J\subset[r]}T_{r, \ell, J}({\bf h}, {\bf k}, {\bf y}, A),
\end{align}
where
\begin{align}
T_{r, \ell, J}({\bf h}, {\bf k}, {\bf y}, A)&:=(-1)^{\wt({\bf h}_{\bar{J}})+\wt({\bf k}_{\bar{I}})+r+|I|}\\
&\quad\times\sum_{\substack{m_j=1\\j\in\bar{J}}}^\infty\prod_{j\in\bar{J}}\frac{e(-m_jy_j)}{m_j^{h_j}}\prod_{i\in\bar{I}}\frac{1}{(\sum_{j\in\bar{J}}a_{ij}m_j)^{k_i}}D({\bf h}_J, {\bf k}_I, {\bf y}_J; \Lambda)\prod_{\substack{j\in J\\i\in I}}\frac{1}{h_j!k_i!}.
\end{align}
\end{theorem}
We can obtain the following corollary by noting $e^{\I\theta}+e^{-\I\theta}=2{\rm Re}(e^{\I\theta})$ for $\theta\in\R$.
\begin{corollary}\label{corollary:1}
We assume that ${\bf h}, {\bf k}, {\bf y}$ and $A$ satisfy the same condition as in \cref{theorem:main theorem parity}. When $\wt({\bf h})+\wt({\bf k})$ and $r$ are of different parity, we have
\begin{align}
{\rm Re}(\zeta_{r, \ell}({\bf h}, {\bf k}, {\bf y},A))=\frac{1}{2}\sum_{\emptyset\neq J\subset[r]}T_{r, \ell, J}({\bf h}, {\bf k}, {\bf y}, A).
\end{align}
\end{corollary}
This means that the real part of $r$-ple series $\zeta_{r, \ell}({\bf h}, {\bf k}, {\bf y}, A)$ can be expressed in terms of $s$-ple series with $s<r$ when $\wt({\bf h})+\wt({\bf k})$ and $r$ are of different parity. Hence, this corollary can be regarded as the parity result for multiple Dirichlet series. Moreover, we can give an expression of ${\rm Re}(\zeta_{r, \ell}({\bf h}, {\bf k}, {\bf y}, A))$ explicitly in terms of $s$-ple series by calculating $D({\bf h}_J, {\bf k}_I, {\bf y}_J; \Lambda)$.
\section{An example}
In this section, we check that the parity result for Mordell--Tornheim type of multiple zeta values is deduced from \cref{corollary:1}.\par
Put $\ell=1, y_1=\dots=y_r=0$ and  
\begin{align}
A=(\underbrace{1, \dots, 1}_{r}).
\end{align}
We can see that $A$ and any positive integers $h_1, \dots, h_r$ and  $k_1$ satisfy the condition of absolute convergence in \cref{theorem:main theorem parity} by the result of Matsumoto and Tsumura~\cite[Lemma 4.2]{MT08}. From now on, we evaluate $T_{r, 1, J}({\bf h}, k_1, {\bf 0}, A)$. In this case, we fix $\rho=(1, 2, \dots, |J|)$ for $J\subset[r]$ and we put
\begin{align}
f_j&=((\delta_{jk})_{k\in J}, 0)\quad(j\in J),\\
f_{r+1}&=((1)_{j\in J},-\sum_{j\in\bar{J}}m_j),\\
\Lambda&=\left\{f_j, f_{r+1}\mid j\in J\right\}.
\end{align}
Then $\mathscr{B}(\Lambda)$ can be divided into $\{B\mid f_{r+1}\not\in B\}$ and $\{B\mid f_{r+1}\in B\}$. The former set contains only one element that is $B=\{f_j\mid j\in J\}$. Hence we have $\vec{f}_j^B=\vec{f}_j$ and $\{{\bf 0}\}_{B, f_j}=0$ for any $j\in J$ since $\ip<\rho, \vec{f}_j^B>>0$.
For $B$ which is an element of the latter set, $B$ corresponds to an element of $J$, we name that element $i$ since there exists only one element of $\Lambda$ such that $f_i\not\in\Lambda\sm B$. Then we can see that
\begin{gather}
\vec{f}_j^B=\vec{f}_j-\vec{f}_i~(j\in J\sm \{i\}), \vec{f}_{r+1}^B=\vec{f}_i,\\
\ip<\vec{f}_i, \vec{f}_j^B>=\begin{cases}\ip<\vec{f}_i, \vec{f}_j-\vec{f}_i>=-1&{\rm if}~j\in J\sm \{i\},\\1&{\rm if}~j=r+1,\end{cases}
\end{gather}
and
\begin{align}
\{{\bf 0}\}_{B, f_j}=\begin{cases}0&{\rm if}~i<j\in(J\cup\{r+1\})\sm \{i\},\\1&{\rm if}~i>j\in(J\cup\{r+1\})\sm \{i\}.\end{cases}
\end{align}
Moreover, $\Z^{|J|}/\ip<\vec{B}>=\{{\bf 0}\}$ for any $B\in\mathscr{B}(\Lambda)$. Therefore we have\par
\begin{align}
G({\bf t}_\Lambda, {\bf 0}; \Lambda)&=\frac{-t_{r+1}}{\dot{f}_{r+1}-\sum_{j\in J}\dot{f}_j\ip<\vec{f}_{r+1}, \vec{f}_j^B>-(t_{r+1}-\sum_{j\in J}t_j\ip<\vec{f}_{r+1}, \vec{f}_j^B>)}\prod_{j\in J}\frac{2\pi\I t_j}{e(t_j)-1}\\
&\quad+\sum_{i\in J}\frac{-t_i}{\dot{f}_i-\sum_{j\in (J\cup\{r+1\})\sm \{i\}}\dot{f}_j\ip<\vec{f}_i, \vec{f}_j^B>-(t_i-\sum_{j\in (J\cup\{r+1\})\sm \{i\}}t_j\ip<\vec{f}_i, \vec{f}_j^B>)}\\
&\quad\times\prod_{\substack{j\in (J\cup\{r+1\})\sm \{i\}\\i<j}}\frac{2\pi\I t_j}{e(t_j)-1}\prod_{\substack{j\in (J\cup\{r+1\})\sm \{i\}\\i>j}}\frac{2\pi\I t_j e(t_j)}{e(t_j)-1}\\
&=\frac{1}{2\pi\I}\frac{1}{\sum_{j\in\bar{J}}m_j-(\sum_{j\in J}t_j-t_{r+1})}\prod_{j\in J\cup\{r+1\}}\frac{2\pi\I t_j}{e(t_j)-1}\\
&\quad\times\left\{e(t_{r+1})-1-\sum_{i\in J}(e(t_i)-1)\prod_{\substack{j\in J\sm\{i\}\\i>j}}e(t_j)\right\}.
\end{align}
Noting that
\begin{align}
\sum_{i\in J}(e(t_i)-1)\prod_{\substack{j\in J\sm\{i\}\\i>j}}e(t_j)=\sum_{i\in J}\left(\prod_{\substack{j\in J\\i\ge j}}e(t_j)-\prod_{\substack{j\in J\sm\{i\}\\i>j}}e(t_j)\right)
=e(\sum_{j\in J}t_j)-1,
\end{align}
we have
\begin{align}
G({\bf t}_\Lambda, {\bf 0}; \Lambda)=\frac{-e(t_{r+1})}{2\pi\I}\frac{e(\sum_{j\in J}t_j-t_{r+1})-1}{\sum_{j\in\bar{J}}m_j-(\sum_{j\in J}t_j-t_{r+1})}\prod_{j\in J\cup\{r+1\}}\frac{2\pi\I t_j}{e(t_j)-1}.
\end{align}
When $J\neq[r]$, we can see that $D({\bf h}_{J}, k_1, {\bf 0}; \Lambda)$ is a polynomial of powers of $\pi\I$ and powers of $(\sum_{j\in\bar{J}}m_j)^{-1}$. When $J=[r]$, we can see that $D({\bf h}_{J}, k_1, {\bf 0}; \Lambda)$ is a polynomial of powers of $\pi\I$. By these observations and by taking the real part, we obtain the parity result for Mordell--Tornheim type of multiple zeta values. 
\begin{remark}
Similarly we can deduce parity results for zeta values of root systems from \cref{corollary:1}.
\end{remark}
\section{Proof of \cref{theorem:main theorem parity}}
Now we prove \cref{theorem:main theorem parity}. Let $M$ be a positive integer. To prove \cref{theorem:main theorem parity}, we evaluate the finite sum
\begin{align}
\zeta_{M, r, \ell}({\bf h}, {\bf k}, {\bf y}, A)=\sum_{m_1=1}^M\dots\sum_{m_r=1}^M\prod_{j\in[r]}\frac{e(m_jy_j)}{m_j^{h_j}}\prod_{i\in[\ell]}\frac{1}{(a_{i1}m_1+\dots+a_{ir}m_r)^{k_i}}.
\end{align}
and then we take the limit $M\to\infty$.\par
From the equation
\begin{align}
\int_0^1e(nx)\,dx=\begin{cases}1&, n=0,\\0&, n\in\Z_{\neq0},\end{cases}
\end{align}
we have
\begin{align}\label{eq:1/c^k}
\frac{1}{c^k}=\suma{n=-N}{N}\frac{1}{n^k}\int_0^1e((c-n)x)\,dx
\end{align}
for a non zero integer $c$ and a positive integer $N$ with $|c|\le N$, where $\sum^*$ means that the dummy variable skips zero. Using \cref{eq:1/c^k} with $c=a_{i1}m_1+\dots+a_{ir}m_r$ and $k=k_i$ for $i\in[\ell]$, we have
\begin{align}
\zeta_{M, r, \ell}({\bf h}, {\bf k}, {\bf y}, A)&=\sum_{m_1=1}^M\dots\sum_{m_r=1}^M \prod_{j\in[r]}\frac{e(m_jy_j)}{m_j^{h_j}}\prod_{i\in[\ell]}\suma{n_i=-aM}{aM}\frac{1}{n_i^{k_i}}\int_0^1e((\sum_{j\in[r]}a_{ij}m_j-n_i)x_i)\,dx_i\\
&=\int_0^1\dots\int_0^1\prod_{j\in[r]}\sum_{m_j=1}^M \frac{e(m_j(y_j+\sum_{i\in[\ell]}a_{ij}x_i))}{m_j^{h_j}}\prod_{i\in[\ell]}\suma{n_i=-aM}{aM}\frac{e(-n_ix_i)}{n_i^{k_i}}\,dx_i,
\end{align}
where $a:=\max\{a_{i1}+\dots+a_{ir}\mid i\in[\ell]\}$. Here, noting the deformation
\begin{align}
\sum_{m=1}^M\frac{e(my)}{m^h}=\suma{m=-M}{M}\frac{e(my)}{m^h}+(-1)^{h+1}\sum_{m=1}^{M}\frac{e(-my)}{m^h},
\end{align}
and expanding the product on $j$, we obtain
\begin{align}
\zeta_{M, r, \ell}({\bf h}, {\bf k}, {\bf y}, A)=\sum_{J\subset[r]}T_{M, r, \ell, J}({\bf h}, {\bf k}, {\bf y}, A),
\end{align}
where $J$ runs over all subsets of $[r]$ and $J$ determines the choice of $j\in[r]$ which $m_j$ runs from $-M$ to $M$ without zero. Moreover, 
\begin{align}
T_{M, r, \ell, J}({\bf h}, {\bf k}, {\bf y}, A)
&:=(-1)^{\wt({\bf h}_{\bar{J}})+|\bar{J}|}\int_0^1\dots\int_0^1\prod_{j\in\bar{J}}\sum_{m_j=1}^M\frac{e(-m_j(y_j+\sum_{i\in[\ell]}a_{ij}x_i))}{m_j^{h_j}}\\
&\quad\times\prod_{j\in J}\suma{m_j=-M}{M}\frac{e(m_j(y_j+\sum_{i\in[\ell]}a_{ij}x_i))}{m_j^{h_j}}\prod_{i\in[\ell]}\suma{n_i=-aM}{aM}\frac{e(-n_ix_i)}{n_i^{k_i}}\,dx_i\\
&=(-1)^{\wt({\bf h}_{\bar{J}})+|\bar{J}|}\prod_{j\in\bar{J}}\sum_{m_j=1}^M\frac{e(-m_jy_j)}{m_j^{h_j}}\prod_{j\in J}\suma{m_j=-M}{M}\frac{e(m_jy_j)}{m_j^{h_j}}\\
&\quad\times\prod_{i\in[\ell]}\suma{n_i=-aM}{aM}\int_0^1e((-\sum_{j\in\bar{J}}a_{ij}m_j+\sum_{j\in J}a_{ij}m_j-n_i)x_i)\,dx_i.
\end{align}
From now on, we evaluate $T_{M, r, \ell, J}({\bf h}, {\bf k}, {\bf y}, A)$ for each $J$. If $\sum_{j\in J}a_{ij}m_j-\sum_{j\in\bar{J}}a_{ij}m_j=0$ for some $i\in[\ell]$, the integral vanishes since $n_i$ skips zero. Hence, $m_j$'s satisfy $\sum_{j\in J}a_{ij}m_j-\sum_{j\in\bar{J}}a_{ij}m_j\neq0$ for all $i\in[\ell]$. From this observation, we obtain
\begin{align}
&T_{M, r, \ell, J}({\bf h}, {\bf k}, {\bf y}, A)\\
&=(-1)^{\wt({\bf h}_{\bar{J}})+|\bar{J}|}\sum_{\substack{m_j=1\\j\in\bar{J}}}^M\prod_{j\in\bar{J}}\frac{e(-m_jy_j)}{m_j^{h_j}}\hspace{-2cm}\sideset{}{^*}\sum_{\substack{m_j=-M\\j\in J\\\sum_{j\in J}a_{ij}m_j-\sum_{j\in\bar{J}}a_{ij}m_j\neq0~({}^\forall i\in[\ell])}}^{M}\hspace{-1.7cm}\prod_{j\in J}\frac{e(m_jy_j)}{m_j^{h_j}}\prod_{i\in[\ell]}\suma{\substack{n_i=-aM\\-\sum_{j\in\bar{J}}a_{ij}m_j+\sum_{j\in J}a_{ij}m_j-n_i=0}}{aM}\frac{1}{n_i^{k_i}}\\
&=(-1)^{\wt({\bf h}_{\bar{J}})+|\bar{J}|}\sum_{\substack{m_j=1\\j\in\bar{J}}}^M\prod_{j\in\bar{J}}\frac{e(-m_jy_j)}{m_j^{h_j}}\hspace{-2cm}\sideset{}{^*}\sum_{\substack{m_j=-M\\j\in J\\\sum_{j\in J}a_{ij}m_j-\sum_{j\in\bar{J}}a_{ij}m_j\neq0~({}^\forall i\in[\ell])}}^{M}\hspace{-1.7cm}\prod_{j\in J}\frac{e(m_jy_j)}{m_j^{h_j}}\prod_{i\in[\ell]}\frac{1}{(\sum_{j\in J}a_{ij}m_j-\sum_{j\in\bar{J}}a_{ij}m_j)^{k_i}}.
\end{align}
Here, we remark that $T_{M, r, \ell, J}({\bf h}, {\bf k}, {\bf y}, A)$ absolutely converges when $M\to\infty$ by the assumption of ${\bf h}, {\bf k}$ and $A$. It is easy to evaluate the case of $J=\emptyset$. We can obtain
\begin{align}
T_{M, r, \ell, \emptyset}({\bf h}, {\bf k}, {\bf y}, A)&=(-1)^{\wt({\bf h})+\wt({\bf k})+r}\zeta_{M, r, \ell}({\bf h}, {\bf k}, -{\bf y}, A)\\
&\xrightarrow{M\to\infty}(-1)^{\wt({\bf h})+\wt({\bf k})+r}\zeta_{r, \ell}({\bf h}, {\bf k}, -{\bf y}, A),\label{eq:zeta to zeta}
\end{align}
where $-{\bf y}=(-y_1, -y_2, \dots, -y_r)$. Hereafter, we assume that $J$ is not the empty set. We see that $T_{M, r, \ell, J}({\bf h}, {\bf k}, {\bf y}, A)$ can be written as
\begin{align}
T_{M, r, \ell, J}({\bf h}, {\bf k}, {\bf y}, A)=(-1)^{\wt({\bf h}_{\bar{J}})+|\bar{J}|}\sum_{\substack{m_j=1\\j\in\bar{J}}}^M\prod_{j\in\bar{J}}\frac{e(-m_jy_j)}{m_j^{h_j}}\hspace{-0.5cm}\sideset{}{^*}\sum_{\substack{m_j=-M\\j\in J\\f_{r+i}({\bf m}_J)\neq0~({}^\forall i\in[\ell])}}^{M}\hspace{-0.5cm}\prod_{j\in J}\frac{e(m_jy_j)}{m_j^{h_j}}\prod_{i\in[\ell]}\frac{1}{f_{r+i}({\bf m}_J)^{k_i}}.
\end{align}
Here, we recall the definition of the set of numbers of non zero row vectors of $A_J$.
\begin{align}
I=I_J:=\{i\in[\ell]\mid a_{ij}\neq0\mbox{ for some }j\in J \},\quad \bar{I}=\bar{I}_J:=[\ell]\sm I_J.
\end{align}
We see that $f_{r+i}({\bf m}_J)=-\sum_{j\in\bar{J}}a_{ij}m_j$ for $i\in\bar{I}$ since $\vec{f}_{r+i}={\bf 0}$ for $i\in\bar{I}$ by the definition of $I$ (see \cref{remark:1}). Moreover, the condition $f_{r+i}({\bf m}_J)\neq0$ obviously holds for $i\in\bar{I}$. Hence, separating the product on $i$ into $I$ and $\bar{I}$, we have
\begin{align}
&T_{M, r, \ell, J}({\bf h}, {\bf k}, {\bf y}, A)\\
&=(-1)^{\wt({\bf h}_{\bar{J}})+|\bar{J}|}\hspace{-0.1cm}\sum_{\substack{m_j=1\\j\in\bar{J}}}^M\prod_{j\in\bar{J}}\frac{e(-m_jy_j)}{m_j^{h_j}}\prod_{i\in\bar{I}}\frac{1}{(-\sum_{j\in\bar{J}}a_{ij}m_j)^{k_i}}\hspace{-0.7cm}\sideset{}{^*}\sum_{\substack{m_j=-M\\j\in J\\f_{r+i}({\bf m}_J)\neq0~({}^\forall i\in I)}}^{M}\hspace{-0.6cm}\prod_{j\in J}\frac{e(m_jy_j)}{m_j^{h_j}}\prod_{i\in I}\frac{1}{f_{r+i}({\bf m}_J)^{k_i}}.
\end{align}
Here, we define
\begin{align}
Z_M({\bf h}_J, {\bf k}_{I}, {\bf y}_J; \Lambda):=\sum_{\begin{subarray}{c}{\bf m}_J\in\Z^{|J|}\\|m_j|\le M~(j\in J)\\f_j({\bf m}_J)\neq0~(j\in J)\\f_{r+i}({\bf m}_J)\neq0~(i\in I)\end{subarray}}e(\ip<{\bf m}_J, {\bf y}_J>)\prod_{j\in J}\frac{1}{f_j({\bf m}_J)^{h_j}}\prod_{i\in I}\frac{1}{f_{r+i}({\bf m}_J)^{k_i}},
\end{align}
then we obtain
\begin{align}
&T_{M, r, \ell, J}({\bf h}, {\bf k}, {\bf y}, A)\\
&=(-1)^{\wt({\bf h}_{\bar{J}})+\wt({\bf k}_{\bar{I}})+|\bar{J}|}\sum_{\substack{m_j=1\\j\in\bar{J}}}^M\prod_{j\in\bar{J}}\frac{e(-m_jy_j)}{m_j^{h_j}}\prod_{i\in\bar{I}}\frac{1}{(\sum_{j\in\bar{J}}a_{ij}m_j)^{k_i}}Z_M({\bf h}_J, {\bf k}_{I}, {\bf y}_J; \Lambda).
\end{align}
By the absolute convergence of $T_{M, r, \ell, J}({\bf h}, {\bf k}, {\bf y}, A)$ and the following result of Komori, Matsumoto and Tsumura \cite[Theorem 2.5]{KMT14CMUSP}
\begin{align}
\lim_{M\to\infty}Z_M({\bf h}_J, {\bf k}_{I}, {\bf y}_J; \Lambda)=(-1)^{|\Lambda|}D({\bf h}_J, {\bf k}_I, {\bf y}_J; \Lambda)\prod_{\substack{j\in J\\i\in I}}\frac{1}{h_j!k_i!},
\end{align}
taking the limit $M\to\infty$, we have
\begin{align}
T_{M, r, \ell, J}({\bf h}, {\bf k}, {\bf y}, A)\to T_{r, \ell, J}({\bf h}, {\bf k}, {\bf y}, A).\label{eq:T to T}
\end{align}
Note that we can use the above result of Komori, Matsumoto and Tsumura for any ${\bf y}_J\in\R^{|J|}$ since the set
\begin{align}
\{f\in\Lambda\mid{\rm rank}\ip<\vec{\Lambda}\sm\{\vec{f}\}>\neq{\rm rank}\ip<\vec{\Lambda}>\}
\end{align}
is empty in our situation. Combining \cref{eq:zeta to zeta,eq:T to T}, we can obtain \cref{theorem:main theorem parity}.

\section*{Acknowledgments}
The author would like to deeply thank Professor Kohji Matsumoto for his moderate encouragement and reading the manuscripts. The author would like to express his sincere gratitude to Professor Yasushi Komori for his valuable advice and comments. The author would like to thank Professor Hirofumi Tsumura for giving him useful comments on a reference. Also, the author would like to thank Dr.~Masataka Ono for his pointing out mistakes in the original version of the article..

\end{document}